\documentclass{article}

\usepackage{cite}

\usepackage{amsmath,amssymb,amsfonts}
\usepackage{algorithmic}
\usepackage{graphicx}
\usepackage{textcomp}
\usepackage{xcolor}
\usepackage{derivative}
\usepackage{url}
\usepackage{caption}
\usepackage{array}
\newcolumntype{P}[1]{>{\centering\arraybackslash}p{#1}}
\newcolumntype{M}[1]{>{\centering\arraybackslash}m{#1}}

\renewcommand{\div}{\nabla \cdot}
\newcommand{\dx}{\mathrm{d}x}
\newcommand{\dt}{\mathrm{d}t}
\newcommand{\ds}{\mathrm{d}s}
\newcommand{\OT}{{\Omega_T}}
\newcommand{\OFEM}{{\Omega_{\mathrm{FEM}}}}
\newcommand{\OFDM}{{\Omega_{\mathrm{FDM}}}}

\newcommand{\ptt}{\partial_{tt}}
\newcommand{\pt}{\partial_{t}}

\newcommand{\pn}{\partial_{n}}

\newcommand{\Eobs}{E_{\mathrm{obs}}} 
\newcommand{\indicator}{\ensuremath{\delta_{\mathrm{obs}}}} 
\newcommand{\z}{\ensuremath{z_{\delta}}} 

\def\BibTeX{{\rm B\kern-.05em{\sc i\kern-.025em b}\kern-.08em
    T\kern-.1667em\lower.7ex\hbox{E}\kern-.125emX}}
\begin{document}

\title{Reconstructing the dielectric properties of melanoma in 3D using real-life melanoma model}


\author{Georg Kyhn,\\
\footnotesize\textit{Department of Mathematical Sciences} \\
\footnotesize\textit{Chalmers University of Technology}\\
\footnotesize\textit{and University of Gothenburg},\\
\footnotesize Gothenburg, Sweden, \\
\footnotesize kyhngeorg@gmail.com
\and
Eric Lindström, \\
\footnotesize\textit{Department of Mathematical Sciences} \\
\footnotesize\textit{Chalmers University of Technology}\\
\footnotesize\textit{and University of Gothenburg},\\
\footnotesize Gothenburg, Sweden, \\
\footnotesize erilinds@chalmers.se
\and
Larisa Beilina, \\
\footnotesize\textit{Department of Mathematical Sciences} \\
\footnotesize\textit{Chalmers University of Technology}\\
\footnotesize\textit{and University of Gothenburg},\\
\footnotesize Gothenburg, Sweden, \\
\footnotesize larisa.beilina@chalmers.se
}

\date{\footnotesize \today}

\maketitle

\begin{abstract}
\footnotesize
  The paper presents performance of the adaptive domain decomposition
  finite element/finite difference method
  for reconstruction of the
  dielectric permittivity and conductivity functions for 3D real-life
  melanoma model using measurements of the backscattered electric field
  at the boundary of the investigated domain. We present several
  gradient-based reconstruction algorithms which use optimization
  approach to find stationary point of the Lagrangian. Our
  computational tests show qualitative and quantitative reconstruction
  of dielectric permittivity and conductivity functions using
  realistic model of malign melanoma at 6 GHz in 3D.

\end{abstract}

\begin{keywords}
  Maxwell’s equations, finite element mesh, adaptive finite element method, finite difference method, coefficient inverse problems, microwave imaging, malign melanoma, dielectric properties of skin
\end{keywords}

\footnotesize{\noindent{\it  \textbf{MSC codes}}: 65J22; 65K10; 65M32; 65M55; 65M60; 65M70}

\graphicspath{
  {/chalmers/groups/waves24/PAPERS/ICEAA2025/Paper2_RecMelanoma/FIGURES/}
    {FIGURES/}
    {pics/}}

\section{Introduction}

This works presents performance of the adaptive domain decomposition
finite element/finite difference method (ADDFE/FDM) for the problem
of determination of the spatially distributed relative dielectric
permittivity and conductivity functions using backscattered data of
the electric field collected at the boundary of the investigated
domain. In scientific community such problems are called Coefficient
Inverse Problems (CIPs) and usually are tackled via minimization of a
least-squares residual functional, see details in
\cite{BakKok},\cite{BK}, \cite{Ghavent}, \cite{T}, \cite{itojin}, \cite{KSS} and
references therein. Globally convergent methods for solution of
electromagnetic CIPs are developed in \cite{convex1,TBKF2,KBKSNF}.

Our mathematical model is described by the stabilized time-dependent
Maxwell's system for the electric field in conductive media which was
theoretically and numerically studied recently in \cite{BR1,BR2},
\cite{BL2,LB3,LB4}. We solve this system numerically via
ADDFE/FDM developed in \cite{BL1}. In works \cite{BL1,BL2,LB4} ADDFE/FDM was
applied for reconstruction of dielectric properties of anatomically
realistic breast phantom at 6 GHz taken from database of
\cite{wisconsin} using adaptive conjugate gradient algorithm (ACGA).

One of the most important application of
proposed algorithms of this work is microwave medical
imaging \cite{ref1}, and particularly, quantitative and qualitative imaging for malign melanoma (MM) detection \cite{noninvasive}.
MM is the
deadliest of skin cancers
though it represents only the 1\% of all skin cancer cases \cite{IEEE2024,ref1}.
The prognosis of MM is related to the primary tumors invasion depth in
the skin \cite{skindepth,ICEAA2024_BEN,IEEE2024}. We note that the
microwave medical imaging
is non-invasive imaging modality compared, for example, to
the X-ray imaging \cite{noninvasive,Pastorino}. It is an attractive addition to the
existing medical imaging technologies like ultrasound
\cite{gonch1,gonch2}, X-ray and MRI imaging \cite{tomography}. We
note that X-ray and MRI are not used in diagnosis of primary skin
cancers, and thus, microwave imaging can be an attractive
non-invasive technology which can be used directly on the skin \cite{IEEE2024}.
Development of this technology is an important compliment to the existing technologies
for diagnosis of MM.

In \cite{IEEE2024,ieee1,wisconsin} were reported different
malign-to-normal tissues contrasts, revealing that malign
tumors have higher relative
permittivity values, than normal tissues, at frequencies less than 10.016 GHz.
It is a big challenge to accurately estimate the relative
permittivity of the internal structure of the skin with MM using the information from the
backscattered electromagnetic waves
collected at several detectors located close to the skin or on the skin with MM.
In the current work we have modelled MM with realistic dielectric properties at 6 GHz taken from \cite{IEEE2024}.

The numerical tests of the current work show that the proposed ADDFE/FDM
method in combination with ACGA reconstruction algorithm can
efficiently and accurately reconstruct the dielectric properties of
malign melanoma at 6 GHz using backscattered noisy data of the
electric field.

\section{The mathematical model}

In this work, the spatial domain of interest $\Omega$ is a bounded, convex and a subset of $\mathbb{R}^3$ with smooth boundary $\Gamma$. 
The time domain of interest is the interval $(0,T)$ for end time $T > 0$. 
Denote the space time domains $\Omega_T := \Omega \times (0,T)$ and $\Gamma_T := \Gamma \times (0,T)$. 
The problem is restricted to linear, isotropic and non-dispersive materials.

In order to utilize the strengths of both the Finite Difference (FD) and the Finite Element (FE) methods, the domain $\Omega$ is decomposed into the subregions $\OFDM$ and $\OFEM$ such that $\Omega = \OFDM \cup \OFEM$ and $\overline{\Omega}_{\mathrm{FEM}} \subset \Omega$ with $\partial \OFEM \subset \OFDM$. 
For details of the domain decomposition, we refer the reader to \cite{BL1}.
Given constants $\overline{\varepsilon}>1$ and $\overline{\sigma}>0$, we assume
\begin{equation}\label{eq:epsilon_sigma_constraint}
    \begin{aligned}
        \varepsilon(x) = 1,& \quad \sigma(x) =  0, && \text{for $x$ in } \OFDM, \\
        \varepsilon(x) \in  [1,\overline{\varepsilon}], & \quad \sigma(x) \in   [0,\overline{\sigma}], && \text{for $x$ in } \Omega \setminus \OFDM,
    \end{aligned}
\end{equation}
holds for the relative electric permittivity $\varepsilon \in C^2(\Omega)$ and the electric conductivity $\sigma \in C^2(\Omega)$.  
From Maxwells equations, utilizing Ohm's, the stabilized problem can be derived. \\  

\noindent
\fbox{
    \begin{minipage}{.97\textwidth}
        \textbf{Forward problem:}
        Find $E$ such that 
        \begin{subequations}\label{eq:forward_prob}
          \begin{align}
              \varepsilon \ptt E + \sigma \pt E - \Delta E - \nabla \div \left( \varepsilon - 1 \right) E 
              &= 0 \ \text{in }\Omega_T, \label{eq:forward_objective_function}\\
              E(\cdot,0) =  f_0, \ \pt E(\cdot,0) 
              &= f_1 \  \text{in } \Omega, \label{eq:forward_initial_cond}\\
              \pn E + \pt E 
              &= 0 \ \text{on } \Gamma_T, \label{eq:forward_absorbing_bc}
          \end{align}
      \end{subequations}
      hold for given $\varepsilon, \sigma, f_0$ and $f_1$. 
    \end{minipage}
}

Above, $E:\OT \to \mathbb{R}^3$ is the electric field and $f_0,f_1:\Omega \to \mathbb{R}$ are  initial conditions.
Note that the Gauss Flux theorem was introduced to the forward problem as a stabilizing term in order to mitigate the possibility of spurious solutions arising due to the use of $P1$-elements in the FEM - see details in \cite{BR1,BR2,BL1,Jiang2}.
The first order absorbing boundary condition \eqref{eq:forward_absorbing_bc}, is included since by \eqref{eq:epsilon_sigma_constraint}, the model equation \eqref{eq:forward_objective_function} becomes the wave equation on $\OFDM$ and $\Gamma \subset \overline{\Omega}_{\mathrm{FDM}}$. 
We can also state the inverse problem. \\

\noindent
\fbox{
    \begin{minipage}{0.97\textwidth}
        \textbf{Inverse problem:}
        Find $\varepsilon(x)$ and $\sigma(x)$ where \eqref{eq:epsilon_sigma_constraint} holds for given bounds $\overline{\varepsilon}$ and $\overline{\sigma}$ such that the following function $\Eobs$ is known at $\Gamma_T$:
        \begin{equation}\label{eq:inverse_problem}
            E \approx \Eobs, \quad \text{on } \Gamma_T.  
        \end{equation}
    \end{minipage}
} \\

Usually, $\Eobs$ are measurements made at the boundary $\Gamma_T$  which incorporate noise as well.

\section{The Lagrangian and optimality conditions}

Due to the fact that the inverse problem is an ill-posed problem, see  \cite{BK,T}, we use an optimization approach and minimize the respective Tikhonov functional, $J$. 
Let
\begin{equation}\label{eq:tikhonov_functional}
    \begin{split}
        J(E,\varepsilon, \sigma ) 
        &:= \frac{1}{2}\int_0^T\int_{\Omega} \left( E - \Eobs \right) ^2 \indicator \z \dx \dt \\
        &+ \frac{\gamma_{\varepsilon}}{2}\int_{\Omega}(\varepsilon - \varepsilon^0)^2\dx 
        + \frac{\gamma_{\sigma}}{2}\int_{\Omega}(\sigma - \sigma^0)^2 \dx. 
    \end{split}
\end{equation}
Here, $\indicator \label{eq:def_indicator}$ is the delta-function equal to one where $\Eobs$ is observed and zero elsewhere, $\z:\mathbb{R} \to \mathbb{R}\label{eq:def_z}$ is a smoothing function, and $\gamma_{\varepsilon}, \gamma_{\sigma} \in \mathbb{R}\label{eq:def_gamma}$ are the regularization parameters for $\varepsilon$ and $\sigma$, respectively. 
Let $H^1(\Omega_T)\label{eq:def_H1}$ be the usual Sobolev space, and define the space $H^1_T(\Omega_T) = \left\{ v \in [H^1(\Omega_T)]^3: v(\cdot,T)= \pt v(\cdot,T) = 0 \right\}$. 
In order to minimize \eqref{eq:tikhonov_functional}, we introduce the Lagrangian $L$ on weak form as 
\begin{equation}\label{eq:lagrangian_weak}
    \begin{split}
        L(E, \lambda, \varepsilon, \sigma)
        &= J(E,\varepsilon,\sigma) 
        -\int_{\Omega} \varepsilon(x) \lambda(x,0) f_1(x) \dx \\
        &- \int_{\Omega} \varepsilon \int_0^T \pt \lambda \pt E \dt \dx 
        + \int_{\Omega} \sigma \int_0^T \lambda \pt E \dt \dx \\
        &+ \int_0^T \int_{\Gamma} \lambda \pt E \ds \dt 
        + \int_0^T\int_{\Omega} \nabla\lambda \nabla E \dx\dt \\
        &+ \int_0^T \int_{\Omega} \div \lambda \div \left( \varepsilon - 1 \right) E \dx\dt.
    \end{split}
\end{equation}
Here, $\lambda \in H^1_T(\Omega_T) \label{eq:def_lambda}$ is the Lagrangian multiplier.
We can define the domain of the Lagrangian as $\mathcal{U} := [H^1(\Omega_T)]^3 \times H_T^1(\Omega_T) \times L^2(\Omega) \times L^2(\Omega)$ such that $L:\mathcal{U} \to \mathbb{R}$. 
Then, to find the minimum of \eqref{eq:tikhonov_functional}, we need to find a stationary point $u \in \mathcal{U}$ of the Lagrangian such that $DL(u)(\tilde{u}) = 0$ for all $\tilde{u} \in \mathcal{U}$. 
Here, $DL(u)$ is the Fréchet derivative of $L$ at $u$.
Since the differentiated Lagrangian can be split into the sum of the partial derivatives, we can consider one term at a time, when finding zeroes. Thus, we obtain:
\begin{align}\label{eq:lagrangian_partial_e}
    \begin{split}
        0=& \ \frac{\partial L}{\partial E}(u)(\tilde{E}) 
        = \int_0^T\int_{\Omega} (E-\Eobs) \tilde{E} \indicator \z \dx \dt \\
        &- \int_\Omega \varepsilon(x) \pt \lambda(x,0)  \tilde{E}(x,0) \dx
        - \int_{\Omega} \varepsilon \int_0^T \pt \lambda \pt \tilde{E} \dt \dx \\
        &- \int_{\Omega} \sigma \int_0^T \pt \lambda \tilde{E} \dt \dx 
       - \int_0^T \int_{\Gamma} \pt \lambda  \tilde{E} \ds \dt \\
        &+ \int_0^T\int_{\Omega} \nabla\lambda \nabla \tilde{E} \dx\dt 
        + \int_0^T \int_{\Omega} \div \lambda \div \left( \varepsilon - 1 \right) \tilde{E} \dx\dt,   
    \end{split}
\end{align}
\begin{align}\label{eq:lagrangian_partial_lambda}
    \begin{split}
        0=& \ \frac{\partial L}{\partial \lambda}(u) (\tilde{\lambda}) 
        = -\int_{\Omega} \varepsilon(x) \tilde{\lambda}(x,0) f_1(x) \dx \\
        &- \int_{\Omega} \varepsilon \int_0^T \pt \tilde{\lambda} \pt E \dt \dx
        + \int_{\Omega} \sigma \int_0^T \tilde{\lambda} \pt E \dt \dx \\
        &+ \int_0^T \int_{\Gamma} \tilde{\lambda} \pt E \ds \dt 
        + \int_0^T\int_{\Omega} \nabla \tilde{\lambda} \nabla E \dx\dt \\
        &+ \int_0^T \int_{\Omega} \div \tilde{\lambda} \div \left( \varepsilon - 1 \right) E \dx\dt,  
    \end{split}
\end{align}
\begin{align}\label{eq:lagrangian_partial_epsilon}
    \begin{split}
       0=\ &\frac{\partial L}{\partial \varepsilon}(u) (\tilde{\varepsilon}) 
        = \gamma_\varepsilon\int_\Omega (\varepsilon - \varepsilon^0)\tilde{\varepsilon} \dx 
        - \int_\Omega \tilde{\varepsilon}(x) \lambda(x,0) f_1(x) \dx \\
        &- \int_\Omega \tilde{\varepsilon} \int_0^T \pt \lambda \pt E \dt \dx 
        + \int_0^T\int_\Omega \nabla \cdot \lambda \nabla \cdot \tilde{\varepsilon} E \dx \dt, 
    \end{split}
\end{align}
\begin{align}\label{eq:lagrangean_partial_sigma}
    \begin{split}
       0= \frac{\partial L}{\partial \sigma}(u) (\tilde{\sigma}) 
        &= \gamma_\sigma\int_\Omega (\sigma - \sigma^0)\tilde{\sigma} \dx 
        + \int_\Omega \tilde{\sigma} \int_0^T \lambda \pt E \dt \dx. 
    \end{split}
\end{align}

Note that finding a stationary point in the $\tilde{\lambda}$ direction is equivalent to solving the forward problem. 
In the same way, from \eqref{eq:lagrangian_partial_e}, the adjoint problem can be derived. \\

\vspace{1mm}
\noindent
\fbox{
    \begin{minipage}{0.97\textwidth}
        \textbf{Adjoint problem:}
        Given $\varepsilon, \sigma$ and $E$, find $\lambda$ such that 
        \begin{subequations}\label{eq:adjoint_prob}
            \begin{align}
                \begin{split}
                    \varepsilon \ptt \lambda &- \sigma \pt \lambda - \Delta \lambda - \left( \varepsilon - 1 \right)\nabla\div \lambda  \\
                    &= \ -(E-\Eobs)\indicator \z,  \label{eq:adjoint_objective_function}
                \end{split}
                && \text{in }\Omega_T,\\
                \lambda(x,T) 
                &= \pt \lambda(x,T) 
                = 0, 
                && \text{in } \Omega, \label{eq:adjoint_end_cond}\\
                \pn \lambda 
                &= \pt \lambda , 
                && \text{on } \Gamma_T, \label{eq:adjoint_absorbing_bc}
            \end{align}
        \end{subequations}
        where $\Eobs$ are observed values on the boundary. 
    \end{minipage}
}

\section{Algorithms for solution of MCIP}

In our numerical simulations we are using
two different gradient-based algorithms: conjugate gradient algorithm (CGA)
 and  adaptive CGA (ACGA).
ACGA uses adaptive
finite element method in space which significantly improves reconstructions of
$\varepsilon, \sigma$ obtained on the initially non-refined mesh by CGA.

Let us define the following functions which are obtained from optimal conditions
$ \pdv{L}{\varepsilon}(u)(\tilde{\varepsilon}) = 0$ and $\pdv{L}{\sigma}(u)(\tilde{\sigma})= 0$
and which we will use in the conjugate gradient algorithm (CGA):
\begin{align}
   & \begin{aligned}\label{gradeps}
        & g_\varepsilon (x) := [\gamma_\varepsilon (\varepsilon - \varepsilon^0)   - \lambda(x,0) f_1 \\
         &- \int_0^T \pt\lambda \pt E ~ \dt
       +  \int_0^T  \div \lambda \div  E~ \dt](x),
     \end{aligned} \\
   & \begin{aligned}
       g_\sigma (x) & := [\gamma_\sigma (\sigma - \sigma^0)  + \int_0^T \lambda \pt E \, \dt] (x).
     \end{aligned}
\end{align}
In \eqref{gradeps} we have used that $\div ( \bar{\varepsilon} E) = \nabla  \bar{\varepsilon} E + \bar{\varepsilon}\div E \approx \bar{\varepsilon}\div E . $

Let $\varepsilon^m$ and $\sigma^m$ be the coefficients computed on the iteration $m$ of the CGA,
and $E^m := E(\varepsilon^m, \sigma^m)$, $\lambda^m := \lambda (E^m, \varepsilon^m, \sigma^m)$, $g_\varepsilon^m := g_\varepsilon^m(x)$  where $E(\varepsilon^m, \sigma^m)$ and $\lambda (E^m, \varepsilon^m, \sigma^m)$ are the solutions of the forward and adjoint problems, respectively.
We introduce CGA  in a continuous setting but the same algorithm
 applies  in the discrete setting.
The conjugate gradient algorithm is formulated as follows. \\

\vspace{.1cm}\noindent\fbox{\begin{minipage}{.97\textwidth}
  \textbf{Conjugate Gradient Algorithm (CGA): }For iterations $m = 0, \dots,M$ perform following steps.
  \begin{enumerate}
    \setcounter{enumi}{-1}
    \item Initialization: set $m = 0$ and choose initial guesses $\varepsilon^0$, $\sigma^0, \gamma_\varepsilon^0, \gamma_\sigma^0, \alpha_\varepsilon^0, \alpha_\sigma^0$.

    \item Calculate $E^m$, $\lambda^m$, $g^m_\varepsilon$ and $g^m_\sigma$.

    \item Update the dielectric permittivity as:
          \begin{align}
            \varepsilon^{m+1}     & := \varepsilon^m + \alpha_\varepsilon^m d_\varepsilon^m,                    \\
            \quad d_\varepsilon^m & := -g_\varepsilon^{m} + \beta_\varepsilon^m d_\varepsilon^{m-1}, 
            \beta_\varepsilon^m    := \frac{\|g_\varepsilon^{m}\|_\Omega^2}{\|g_\varepsilon^{m-1}\|_\Omega^2},
          \end{align}

    \item Update conductivity as:
          \begin{flalign}
            \sigma^{m+1}     & := \sigma^m + \alpha_\sigma^m d_\sigma^m,     &&                     \\
            \quad d_\sigma^m & := -g_\sigma^{m} + \beta_\sigma^m d_\sigma^{m-1}, 
            \beta_\sigma^m    := \frac{\|g_\sigma^{m}\|_\Omega^2}{\|g_\sigma^{m-1},\|_\Omega^2}, &&
          \end{flalign}
          where $\alpha_\varepsilon$, $\alpha_\sigma$ are chosen step sizes with $d_\varepsilon^0 := - g_\varepsilon^0$, $d_\sigma^0 := - g_\sigma^0$.

    \item Compute new optimal step-sizes as
          \begin{equation}
            \begin{split}
              \alpha_\varepsilon^{m+1} = - \frac{(g_\varepsilon^m, d_\varepsilon^m)}{ \gamma_\varepsilon^m ( d_\varepsilon^m, d_\varepsilon^m )}, \quad
              \alpha_\sigma^{m+1} =- \frac{(g_\sigma^m, d_\sigma^m)}{ \gamma_\sigma^m ( d_\sigma^m, d_\sigma^m )}.
            \end{split}
          \end{equation}

    \item Compute new regularization parameters for any $p \in (0,1)$ via iterative rules of \cite{BakKok} as
          \begin{equation}
            \begin{split}
              \gamma_\varepsilon^{m+1} = \frac{ \gamma_\varepsilon^0}{(m+1)^p},  \quad
              \gamma_\sigma^{m+1} =  \frac{ \gamma_\sigma^0}{(m+1)^p}
            \end{split}
          \end{equation}

    \item Terminate the algorithm if either $\|\varepsilon^{m+1} - \varepsilon^m \| < \eta_\varepsilon^1$ or $\|\sigma^{m+1} - \sigma^m \| < \eta_\sigma^1$ and $\|g_\varepsilon^m\| < \eta_\varepsilon^2$ or $\|g_\sigma^m\| < \eta_\sigma^2$, where $\eta^1_\varepsilon$, $\eta^2_\varepsilon$, $\eta^1_\sigma$ and $\eta^2_\sigma$ are tolerances chosen by the user. Otherwise, set $m: = m+1$ and repeat the algorithm from step 1).
  \end{enumerate}
\end{minipage}
}\\

The main idea of the adaptive local mesh refinement used in this work
is that the finite element mesh $K_h$ should be refined in such elements where
 both functions
$|h {\varepsilon}_h|, |h {\sigma}_h| $ achieve its maximum
values:
\begin{equation}\label{meshref}
  \begin{split}
    |h {\varepsilon}_h|
    \geq \widetilde{\beta}_\varepsilon \max \limits_{K_h} |h {\varepsilon}_h|, \\
    |h {\sigma}_h| \geq \widetilde{\beta}_\sigma \max \limits_{K_h} |h {\sigma}_h|.
  \end{split}
\end{equation}
Here, $ \widetilde{\beta}_\varepsilon \in (0,1) , \widetilde{\beta}_\sigma \in (0,1)$ are
numbers which should
be chosen computationally, and $h = h(x)$ is a piecewise-constant mesh function
representing the local diameter of the elements and is defined as
\begin{equation}\label{meshfunction}
  h |_K = h_K ~~~ \forall K \in K_h.
\end{equation}

The mesh refinements recommendations \eqref{meshref} are based on a posteriori error estimates for the
errors $| \varepsilon - {\varepsilon}_h|$, $| \sigma -
 \sigma_h|$ in the reconstructed functions $\varepsilon, \sigma$,
respectively. The proofs of  a posteriori error estimates
 can be derived using technique of \cite{BL2}
and
is
topic of the ongoing research. We formulate our adaptive conjugate
gradient algorithm as follows below.

\vspace{.1cm}\noindent\fbox{\begin{minipage}{.97\textwidth}
    \textbf{Adaptive Conjugate Gradient Algorithm (ACGA): }For mesh
    refinements $i = 0, \dots, N$ perform the steps below.
    \begin{enumerate} \setcounter{enumi}{-1} \item Choose initial
            spatial mesh $K^0_h$ in $\OFEM$.

      \item Compute $\varepsilon_{h,i}$, $\sigma_{h,i}$ on mesh $K^i_h$ according to the CGA algorithm.

      \item Refine locally such elements in the finite element mesh $K_h^i$  where
            \begin{equation*}
              \begin{split}
                |h {\varepsilon}_{h,i}|
                \geq \widetilde{\beta}_{\varepsilon,i}  \max_{K  \in K_h^i} |h {\varepsilon}_{h,i}|, \\
                |h {\sigma}_{h,i}| \geq \widetilde{\beta}_{\sigma,i}  \max_{K  \in K_h^i}  |h {\sigma}_{h,i}|.
              \end{split}
            \end{equation*}
            Here, $ \widetilde{\beta}_{\varepsilon,i} \in (0,1)  , \widetilde{\beta}_{\sigma,i} \in (0,1)$
            are constants chosen by the user.

      \item Construct the new mesh as $K_h^{i+1}$ and interpolate $\varepsilon_{h,i}$, $\sigma_{h,i}$ as well as measurements $E_{\rm obs}$ onto it.

      \item Terminate the algorithm if either $\|\varepsilon_{h,i+1} - \varepsilon_{h,i} \| < \theta_\varepsilon^1$ or $\|\sigma_{h, i+1} - \sigma_{h,i} \| < \theta_\sigma^1$ and $\|g_\varepsilon^m\| < \theta_\varepsilon^2$ or $\|g_\sigma^m\| < \theta_\sigma^2$, where $\theta^1_\varepsilon$, $\theta^2_\varepsilon$, $\theta^1_\sigma$ and $\theta^2_\sigma$ are tolerances chosen by the user. Otherwise, set $i := i+1$ and repeat the algorithm from step 1).
    \end{enumerate}
  \end{minipage}
} \\

We note that our mesh refinement is performed only in space and not in
time to be able smoothly run the ADDFE/FDM algorithm
for solutions of the forward and adjoint problems - see details of
this algorithm for our model problem in \cite{BL1}. We also note that
numbers $ \widetilde{\beta}_{\varepsilon,i}, \widetilde{\beta}_{\sigma,i}$ in ACGA algorithm should be chosen computationally depending
on the obtained results of reconstruction. If these numbers are
chosen close to zero then almost all elements in the mesh $K_h^i$ will
be refined, and opposite, if $\widetilde{\beta}_{\varepsilon,i} \approx 1, \widetilde{\beta}_{\sigma,i}\approx 1$ then the mesh will be not refined
at all.

\begin{figure}[tbp]
  \begin{center}
    \begin{tabular}{cc}
      \includegraphics[trim = 8.0cm 4.0cm 4.0cm 4.0cm, scale=0.16, clip=]{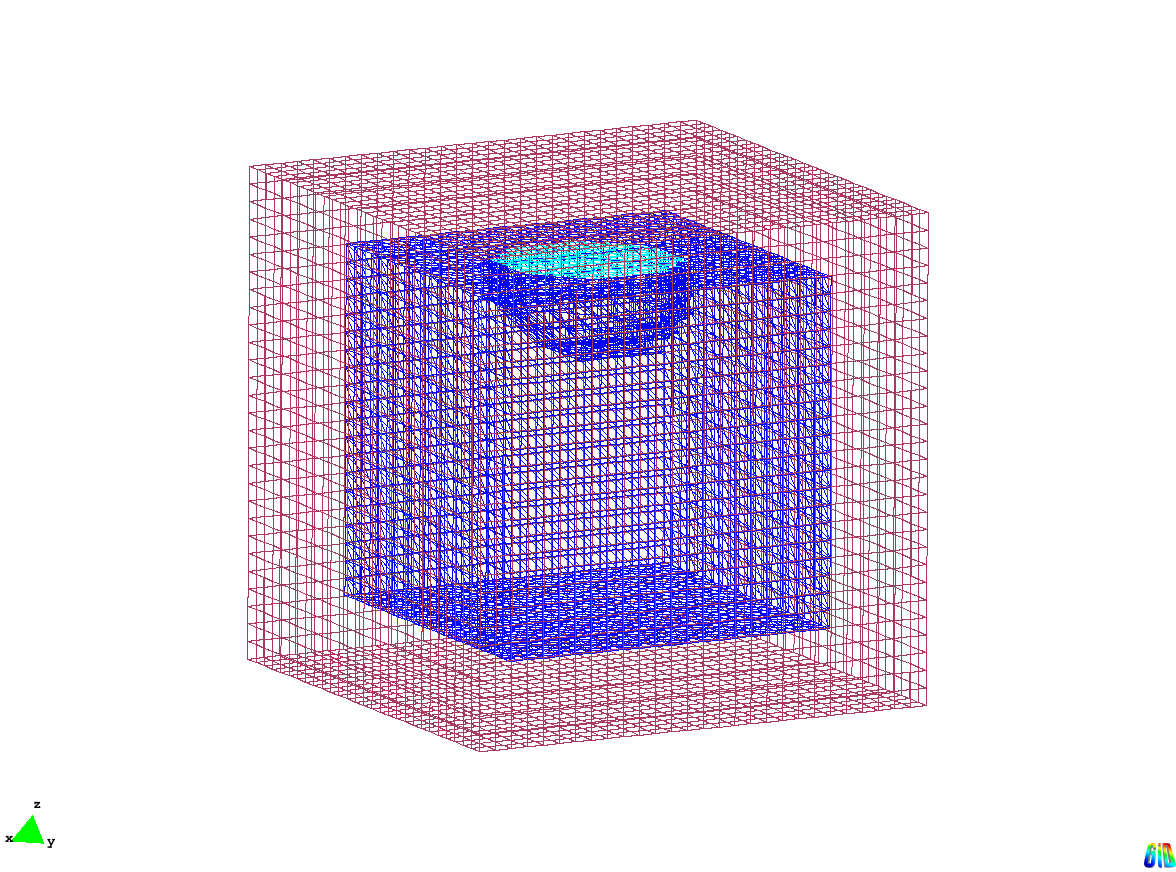} &
      \includegraphics[trim = 8.0cm 3.0cm 4.0cm 3.0cm, scale=0.16, clip=]{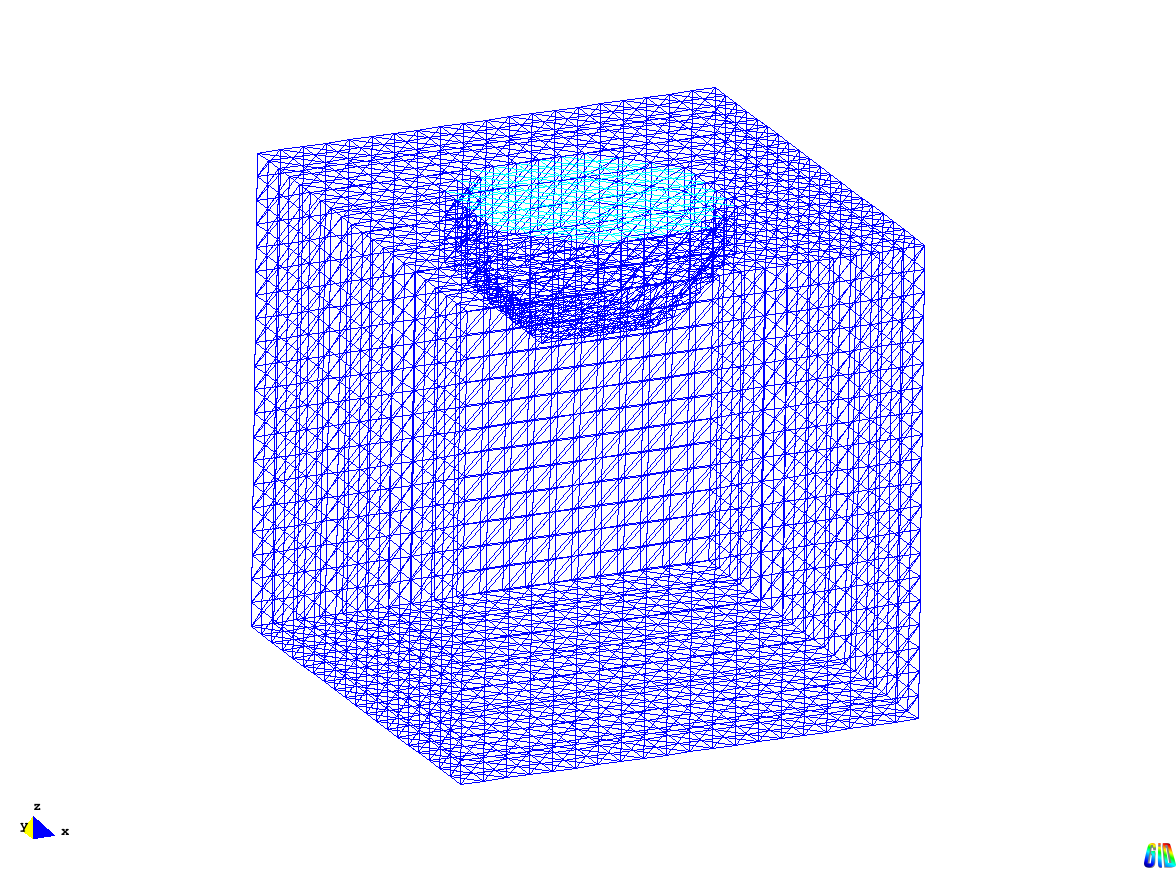}
    \end{tabular}
  \end{center}
  \caption{\small
    \textit{Computational domains used in the domain decomposition ADDFE/FDM.
    Left: the hybrid finite element/finite difference domain $\Omega : = \Omega_{\rm FEM} \cup \Omega_{\rm FDM}$. Right: the finite element domain $\Omega_{\rm FEM}$ with model of malign melanoma.}}
  \label{fig:meshes}
\end{figure}

\section{Numerical results}

This section presents numerical results of the reconstruction of the
relative dielectric permittivity $\varepsilon$ and conductivity $\sigma$ using
CGA and ACGA algorithms. We refer to \cite{BL1,BL2} for details about
numerical implementation of the ADDFE/FDM method which was used for
solutions of the forward and adjoint problems.

\begin{figure}[tbp]
  \begin{center}

    \begin{tabular}{cc}
      \includegraphics[trim = 2.0cm 1.0cm 0.0cm 5.0cm, scale=0.16, clip=]{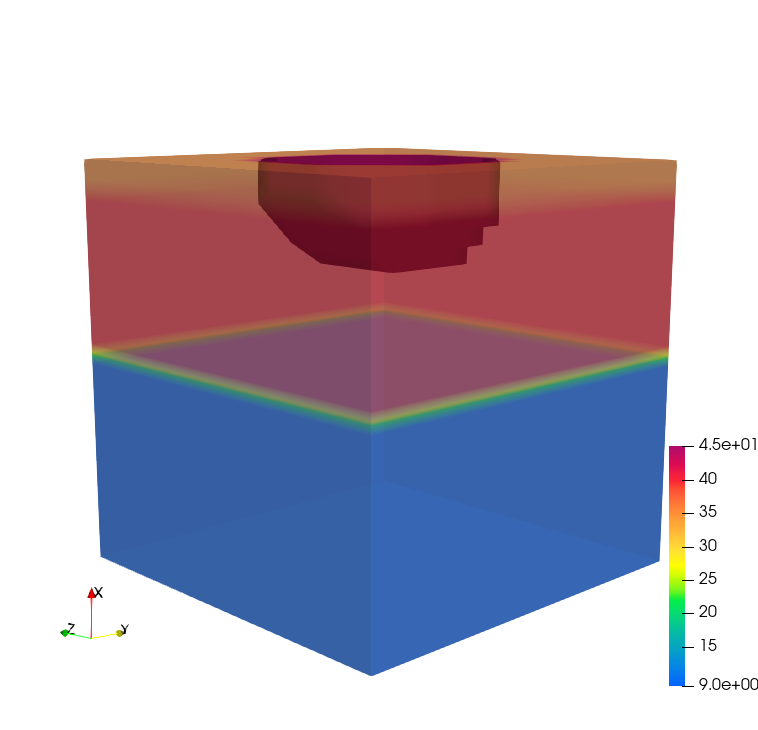} &
      \includegraphics[trim = 2.0cm 1.0cm 0.0cm 5.0cm, scale=0.16, clip=]{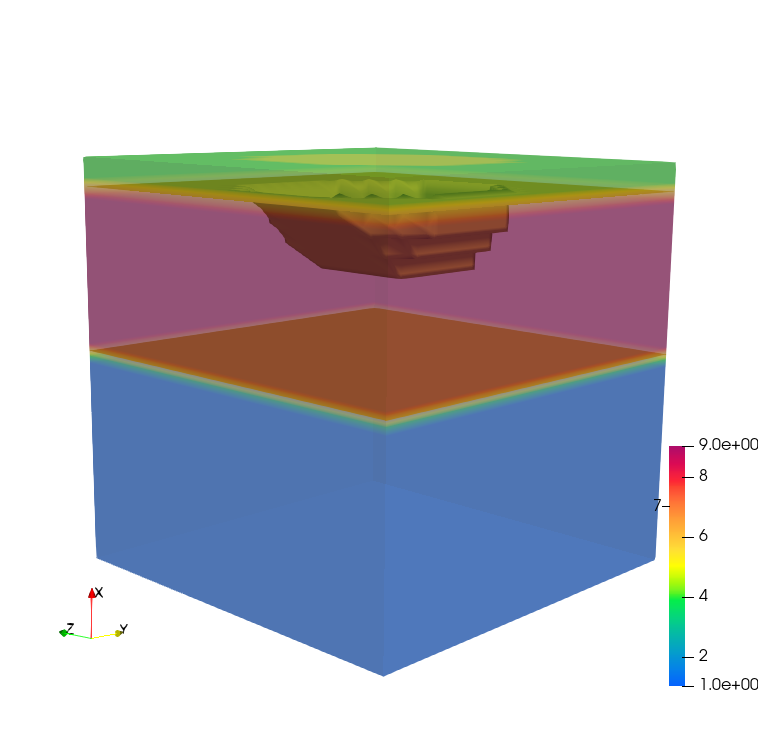}                                    \\
      $\max_\Omega \varepsilon_{r} \approx 45$                                                                   & $\max_\Omega \sigma \approx 9$ \\
    \end{tabular}
  \end{center}
  \vspace{-.4cm}
  \caption{\small
    \textit{Model with realistic dielectric properties of malign melanoma and skin at 6 GHz corresponding to the exact values of
    $\varepsilon_{r}, \sigma$ in the Table \ref{tab:table1}}.
  }
  \label{fig:numex1}
\end{figure}

\begin{figure}[h!]
  \begin{center}
  \vspace{-.2cm}
    \begin{tabular}{cccc}
      \includegraphics[trim = 1.0cm 0.0cm 1.0cm 4.0cm, scale=0.11, clip=]{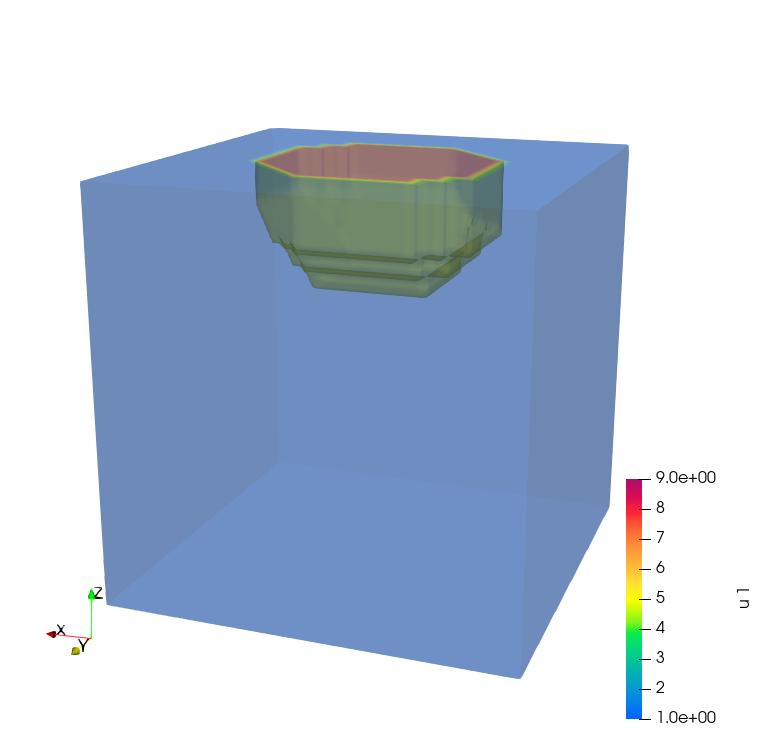} &
      \includegraphics[trim = 1.0cm 0.0cm 1.0cm 4.0cm, scale=0.11, clip=]{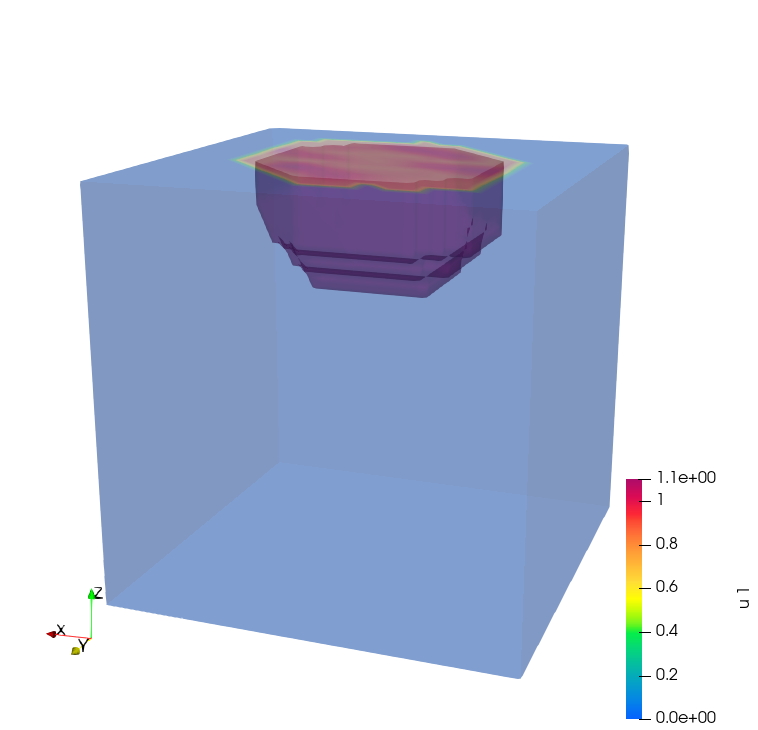} &
      \includegraphics[trim = 10.0cm 4.0cm 11.5cm 4.0cm, scale=0.11, clip=]{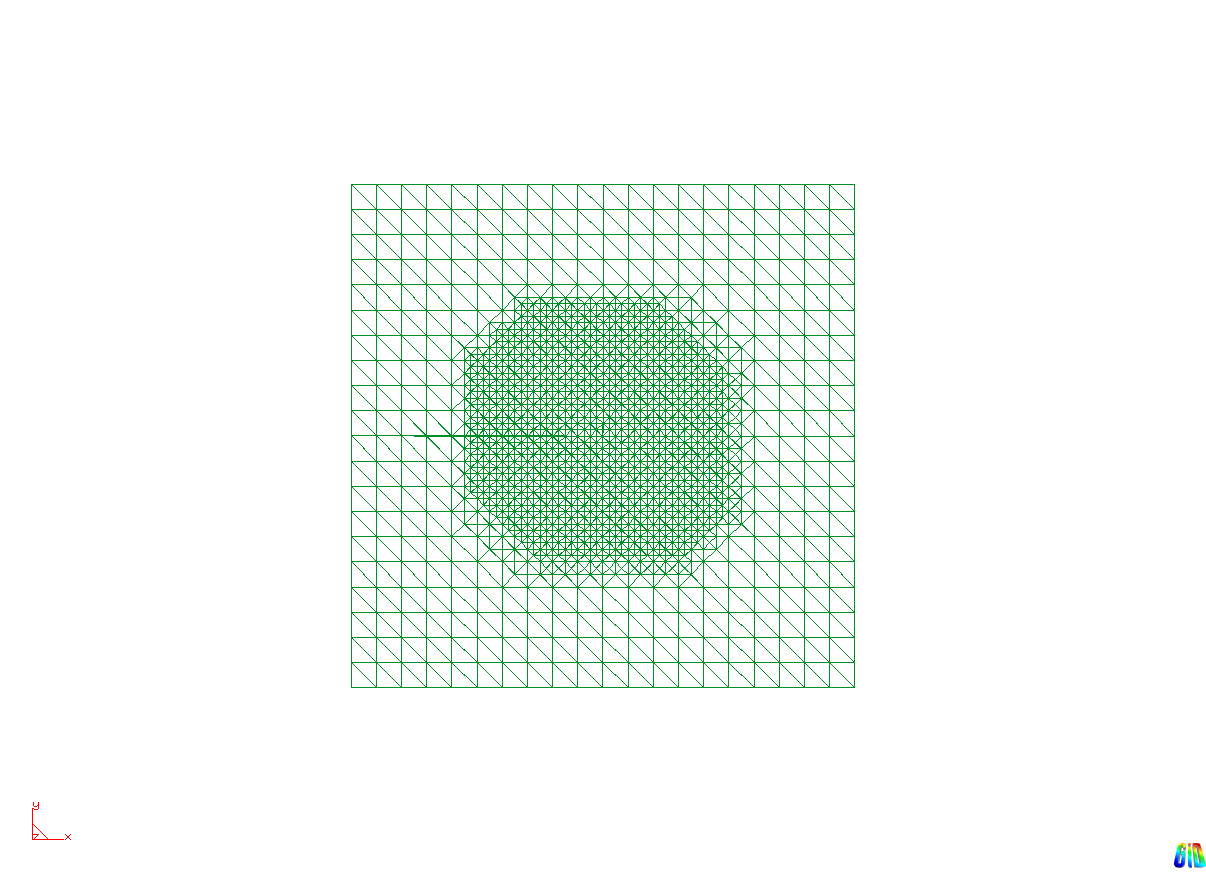}   &
      \includegraphics[trim = 8.0cm 4.0cm 11.5cm 4.0cm, scale=0.11, clip=]{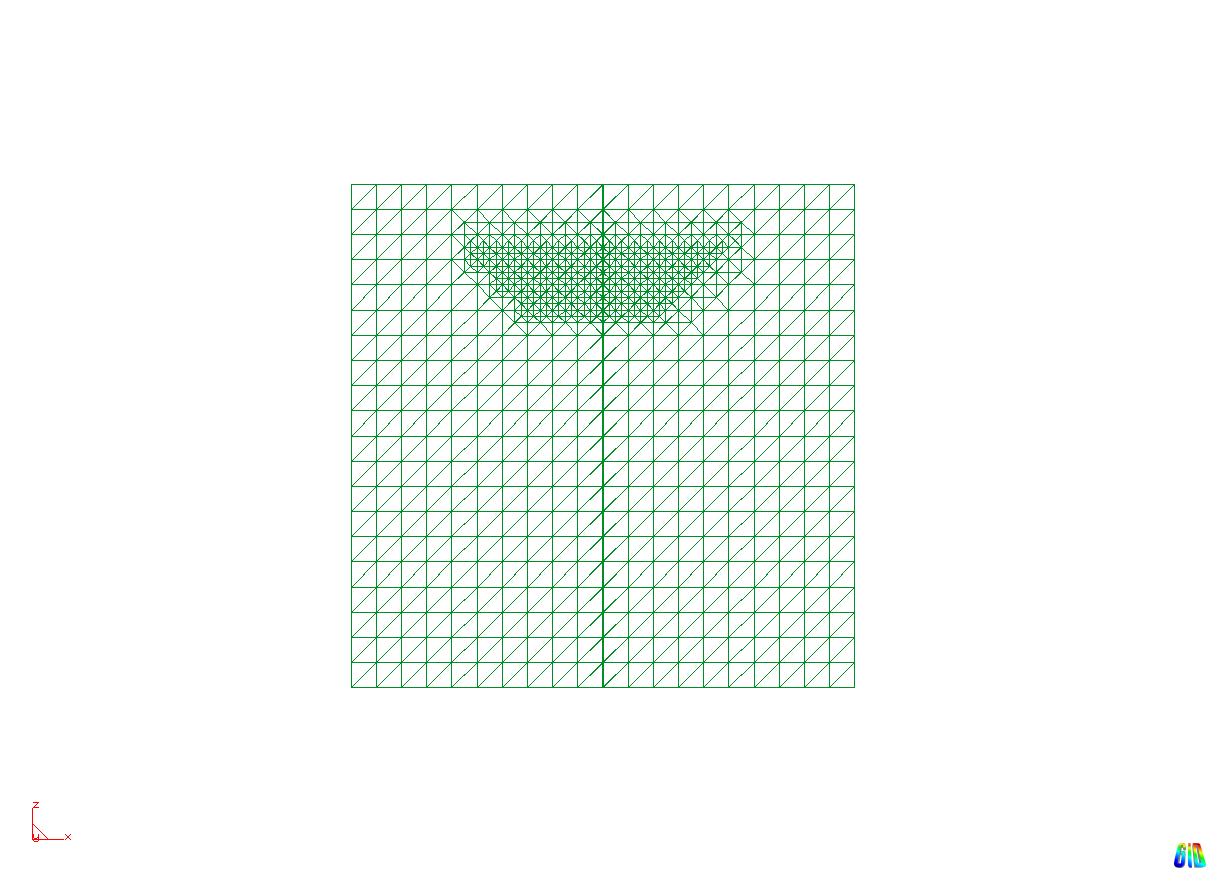}                                                 \\
      \small a)  $  \max_{K_{h_5}}  \varepsilon_{r}/5  \approx 9$  & b)   $ \max_{K_{h_5}} \sigma/5 \approx 1.2$  &c)  $x_1 x_2$-view                                                                                 & d) $x_3 x_1$- view
    \end{tabular}
  \end{center}
  \vspace{-.5cm}
  \caption{\small
    \textit{Dielectric weighted properties at 6GHz taken in our computations: a) $\varepsilon$, b) $\sigma$.
    c), d) Projections of the locally refined mesh ${K_{h_5}}$ taken for generation of data.}
  }
  \label{fig:numex2}
\end{figure}

We use the same model of malign melanoma as was developed in
\cite{ICEAA2024_BEN}. To be able apply ADDFE/FDM for this model we
split the computational domain $\Omega$ into two domains $\Omega_{\rm FEM},
 \Omega_{\rm FDM}$ such that
$\Omega : = \Omega_{\rm FEM} \cup \Omega_{\rm FDM}$ with
$\Omega_{\rm FEM} \subset \Omega$.
The domain $\Omega_{\rm FEM}$
is discretized by
tetrahedral elements and the domain $\Omega_{\rm FDM}$ - by hexahedral elements - see Figure \ref{fig:meshes} and \cite{BL1} for further
details of discretization in the domain decomposition method.
More precisely, we set
the dimensionless computational domain $\Omega$ as
\begin{equation*}
 \Omega = \left\{ x= (x_1,x_2, x_3) \in (-2, 12) \times (-2, 12) \times (-2,
 12) \right\},
\end{equation*}
and the computational domain $\Omega_{\rm FEM}$ as
\begin{equation*}
\Omega_{\rm FEM} =
 \left\{ x= (x_1,x_2, x_3) \in (0, 10) \times (0, 10) \times (0,
 10) \right\}.
\end{equation*}
The domain $\Omega_{\rm FEM}$ corresponds to the 3D model of MM of the size $10 \times 10 \times 10$ mm.
We assign values of $\varepsilon$
and $\sigma$ in $\Omega_{\rm FEM}$ accordingly to the test values of the Table
\ref{tab:table1}, and we set $\varepsilon =1$
and $\sigma = 0$ in $\Omega \setminus \Omega_{\rm FEM}$.
Figure \ref{fig:numex1} presents exact values of the relative dielectric permittivity and conductivity functions at 6 GHz corresponding to the Table 1 for the real skin model with all tissue types (epidermis, dermis, fat). Figure \ref{fig:numex2} shows simplified model of
Figure \ref{fig:numex1} with exact values of the relative dielectric permittivity and conductivity functions which we are reconstructing
in our numerical tests. Dielectric properties shown on this figure
have weighted values of $\varepsilon$ and $\sigma$ which correspond to the test values presented in the Table \ref{tab:table1}.
To proceed further, we decompose the boundary $\partial \Omega$ of the domain $\Omega$ into
three parts as follows: $\partial
 \Omega =\partial _{1} \Omega \cup \partial _{2} \Omega \cup \partial
 _{3} \Omega$. Here, $\partial _{1} \Omega$ and $\partial _{2} \Omega$
are, respectively, top and bottom sides of $\Omega$, and $\partial
 _{3} \Omega$ is the union of left, right, front and back sides of this
domain. The time-dependent observations are simulated at
the backscattered boundary
$\Gamma_1 :=
 \partial_1 \Omega \times (0,T)$ of $\Omega$. We also define $\Gamma_{1,1} := \partial_1
 \Omega \times (0,t_1]$, $\Gamma_{1,2} := \partial_1 \Omega \times
 (t_1,T)$, and $\Gamma_3
 := \partial_3 \Omega \times (0, T)$.
In our computations we use the following stabilized
model problem:
\begin{equation}\label{model1}
  \begin{split}
    \varepsilon \ptt E & +
    \nabla (\nabla \cdot E) - \Delta E  -
    \nabla  (\nabla \cdot ( \varepsilon  E))
    = -  \sigma \pt E                              ~\mbox{in}~ \Omega_T,                   \\
    E(x,0)                & = 0, ~~~\pt E(x,0) = 0 ~\mbox{in}~ \Omega,                     \\
    \pn E                 & = P(t)                 ~\mbox{on}~ \Gamma_{1,1},               \\
    \pn E                 & = - \pt E              ~\mbox{on}~ \Gamma_{1,2} \cup \Gamma_2, \\
    \pn E                 & = 0                    ~\mbox{on}~ \Gamma_3.                   \\
  \end{split}
\end{equation}
We initialize a plane wave $P(t) = (0,P_2,0)(t)$ only for the one component $E_2$ of the electric field
$E=(E_1, E_2, E_3)$ at $\Gamma_{1,1}$
in (\ref{model1})
in time $t=[0,12.0]$ and define it as
\begin{equation}\label{f}
  \begin{split}
    P_2(t) =\left\{
    \begin{array}{ll}
      \sin \left( \omega t \right) ,\qquad & \text{ if }t\in \left( 0,\frac{2\pi }{\omega }
      \right) ,                                                                             \\
      0,                                   & \text{ if } t>\frac{2\pi }{\omega }.
    \end{array}
    \right.
  \end{split}
\end{equation}

\begin{table}[ht]
  \centering
  \small
  \caption{\textit{Tissue types and corresponding realistic
      values of $\varepsilon$ and $\sigma$ (S/m) at 6 Ghz for skin with melanoma
      used in our numerical experiments. }}
  \begin{tabular}{| p{2.3cm} | P{.5cm} | P{.6cm}  | P{1cm} |  P{1cm}  | P{.7cm} | }
    \hline
                     & \multicolumn{2}{c|}{Real values} & \multicolumn{2}{c|}{Test values} &                                                                \\
    \hline
    Tissue type      & $\varepsilon$                    & $\sigma$                         & $\varepsilon\,(\varepsilon/5)$ & $\sigma\,(\sigma/5)$ & Depth  \\
                     &                                  & (S/m)                            &                                & (S/m)                & (mm)   \\
    \hline
    \hline
    Immersion medium & 32                               & 4                                & 5 (1)                          & 0                    & 2      \\
    \hline
    Epidermis        & 35                               & 4                                & 5 (1)                          & 0                    & 1      \\
    \hline
    Dermis           & 40                               & 9                                & 5 (1)                          & 0                    & 3.5    \\
    \hline
    Fat              & 9                                & 1                                & 5 (1)                          & 0                    & 5.5    \\
    \hline
    Tumor stage 1    & 45                               & 6                                & 40 (8)                         & 6 (1.2)              & $< 1$  \\
    \hline
    Tumor stage 2    & 50                               & 6                                & 45 (9)                         & 6 (1.2)              & $ > 1$ \\
    \hline
    Tumor stage 3    & 60                               & 6                                & -                              & -                    & $> 1$  \\
    \hline
  \end{tabular}
  \label{tab:table1}
\end{table}

\begin{figure*}[h!]
  \begin{center}
    \begin{tabular}{ccccc}
      \includegraphics[trim = 2.0cm .2cm 1.0cm 4.0cm, scale=0.11, clip=]{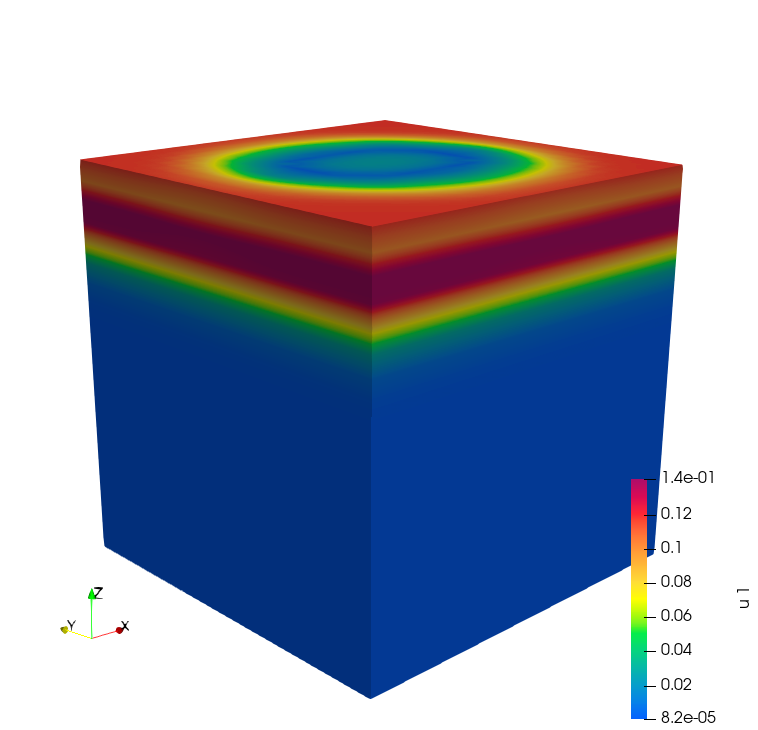} &
      \includegraphics[trim = 2.0cm .2cm 1.0cm 4.0cm, scale=0.11, clip=]{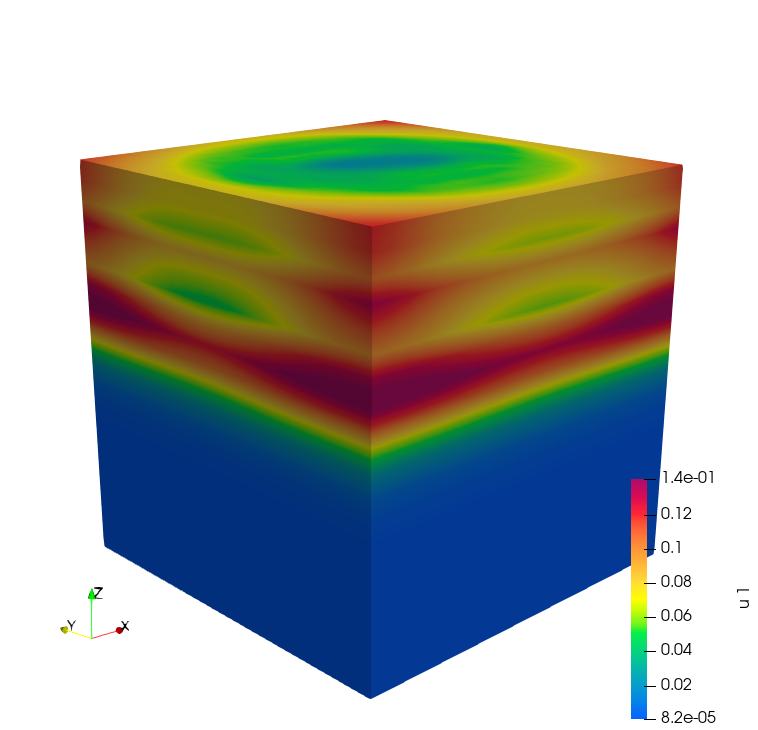} &
      \includegraphics[trim = 2.0cm .2cm 1.0cm 4.0cm, scale=0.11, clip=]{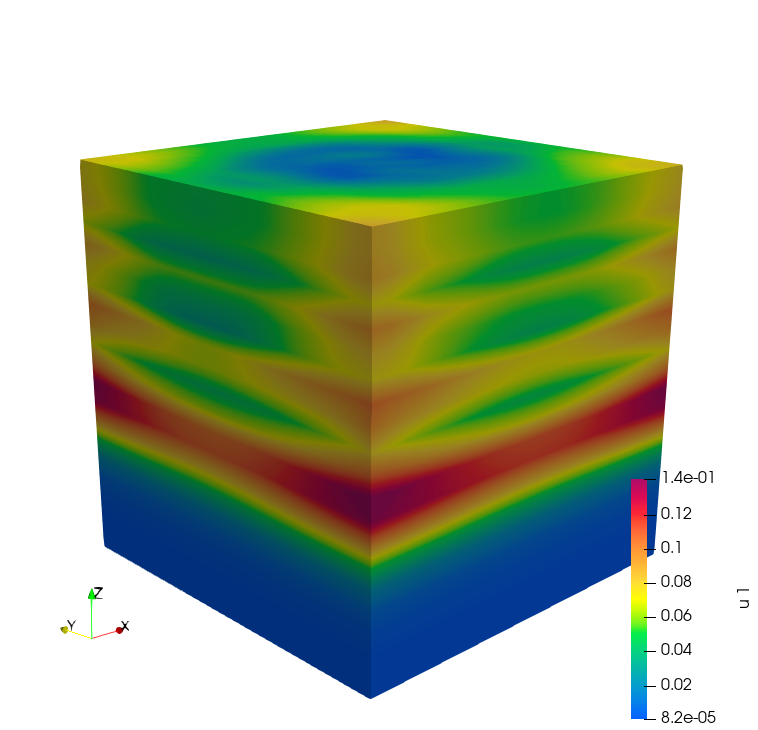} &
      \includegraphics[trim = 2.0cm .2cm 1.0cm 4.0cm, scale=0.11, clip=]{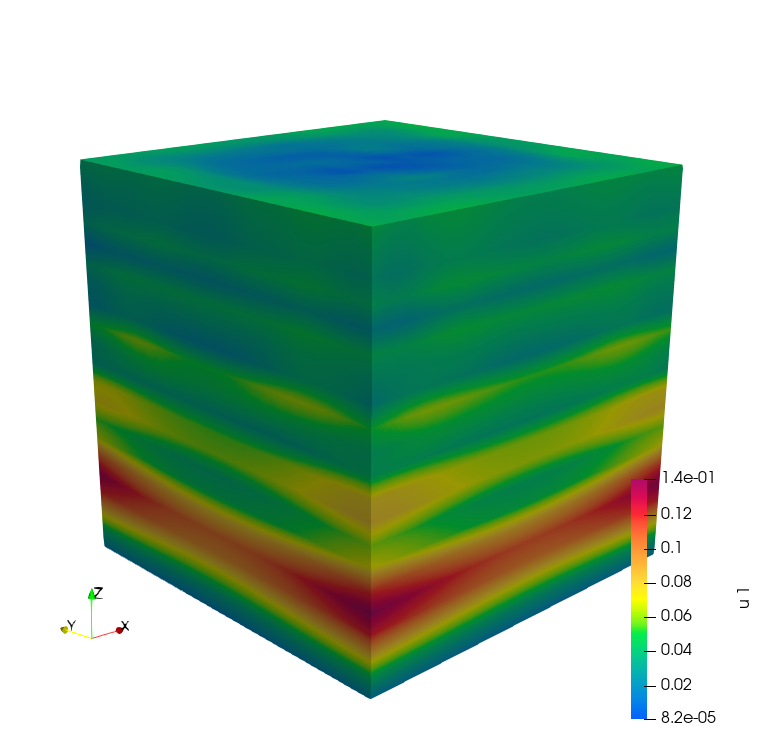}
    \end{tabular}
  \end{center}
  \caption{\small
    \textit{Generation of data. Dynamics of scattered electric field $|E|$ in the finite element domain at different times using DD FE/FD method computed on the geometry shown on Figure \ref{fig:meshes}. We note that in the reconstruction algorithm were used only backscattered data at the top boundary.}
  }
  \label{fig:data}
\end{figure*}

\begin{figure}[h!]
  \begin{center}
    \begin{tabular}{cccc}
      \includegraphics[trim = 3.0cm 0.0cm 1.5cm 4.0cm, scale=0.1, clip=]{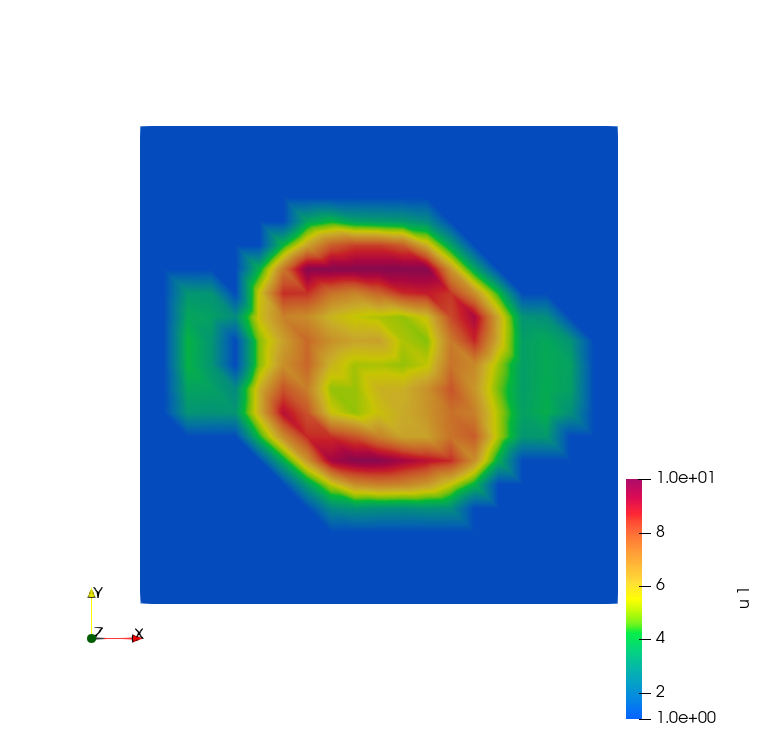} &
      \includegraphics[trim = 3.0cm 0.0cm 1.5cm 4.0cm, scale=0.1, clip=]{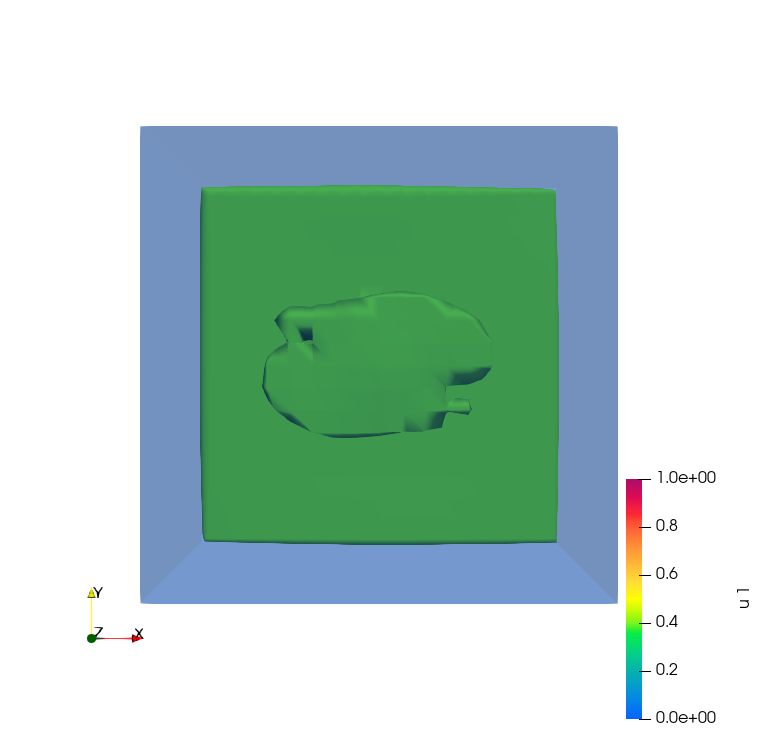} &
      \includegraphics[trim = 0.0cm 0.0cm 1.0cm 0.0cm, scale=0.1, clip=]{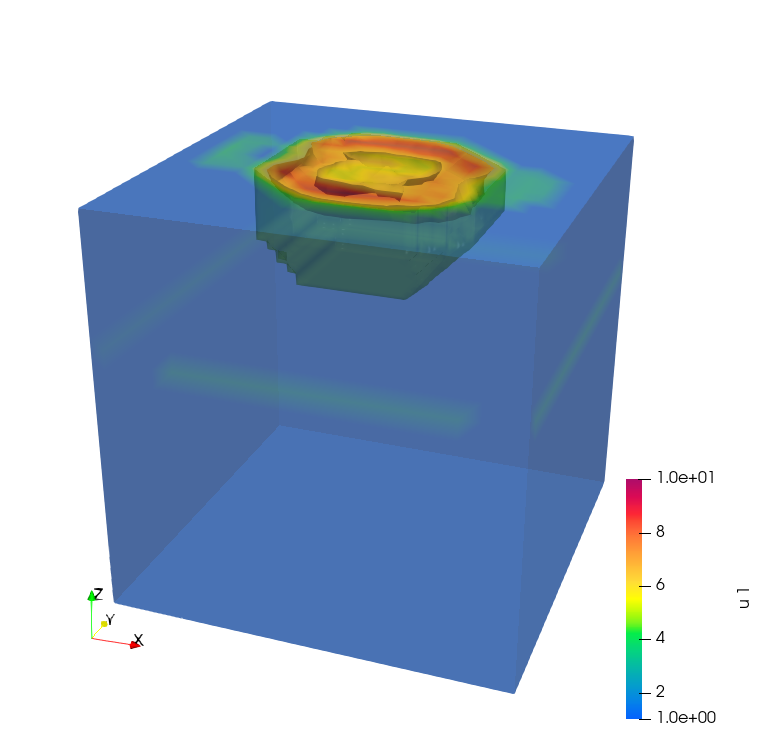}    &
      \includegraphics[trim = 0.0cm 0.0cm 1.0cm 0.0cm, scale=0.1, clip=]{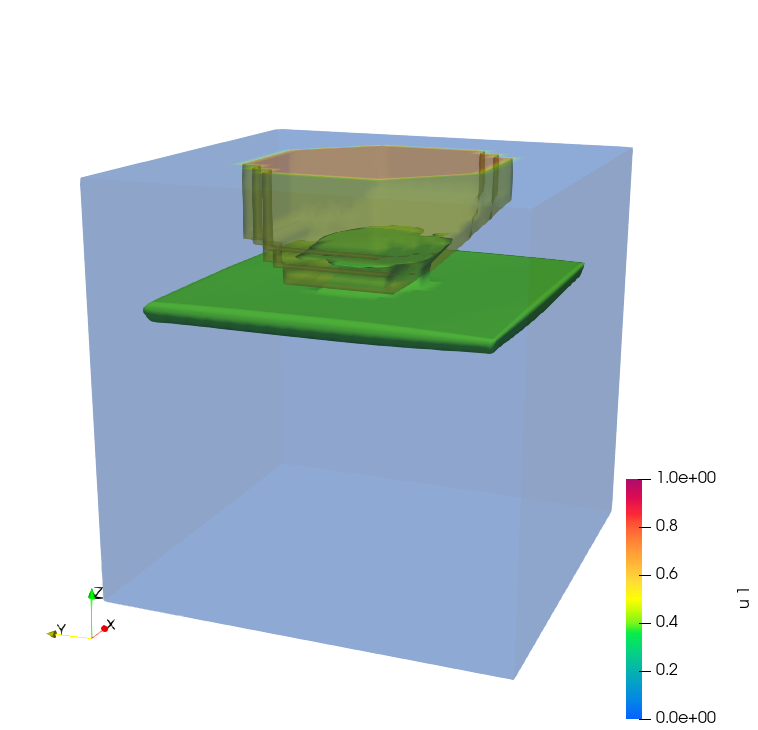}
      \\
      \small a) $\max \varepsilon_{r_h}/5  \approx 9$ & b) $\max \sigma_{h}/5 \approx 0.6$ &c) $\max \varepsilon_{r_h}/5  \approx 9$                                                                       & d) $\max \sigma_{h}/5 \approx 0.6$
    \end{tabular}
  \end{center}
  \vspace{-.4cm}
  \caption{\small \textit{Performance of CGA on the coarse mesh ${K_h}_0$.
    a),c) The weighted reconstruction of $\varepsilon_{r_h}$ (outlined in red and yellow colors) corresponding to
    the malign melanoma at stage 1.
    b),d) The weighted reconstruction of $\sigma_{h}$ (outlined in green color).
    The noise level in the data for electric field is $\delta= 10\%$.}
  }
  \label{fig:CGA}
\end{figure}

\begin{figure}[h!]
  \begin{center}
    \begin{tabular}{cc}
      \includegraphics[trim = 2.7cm 0.2cm 1cm 4.0cm, scale=0.11, clip=]{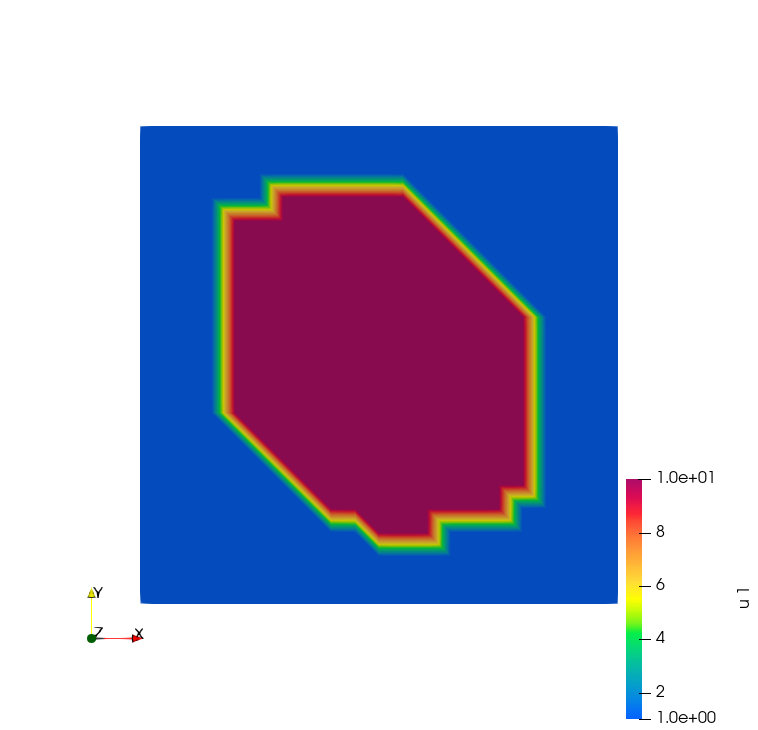} &
      \includegraphics[trim = 2.7cm 0.2cm 1cm 4.0cm, scale=0.11, clip=]{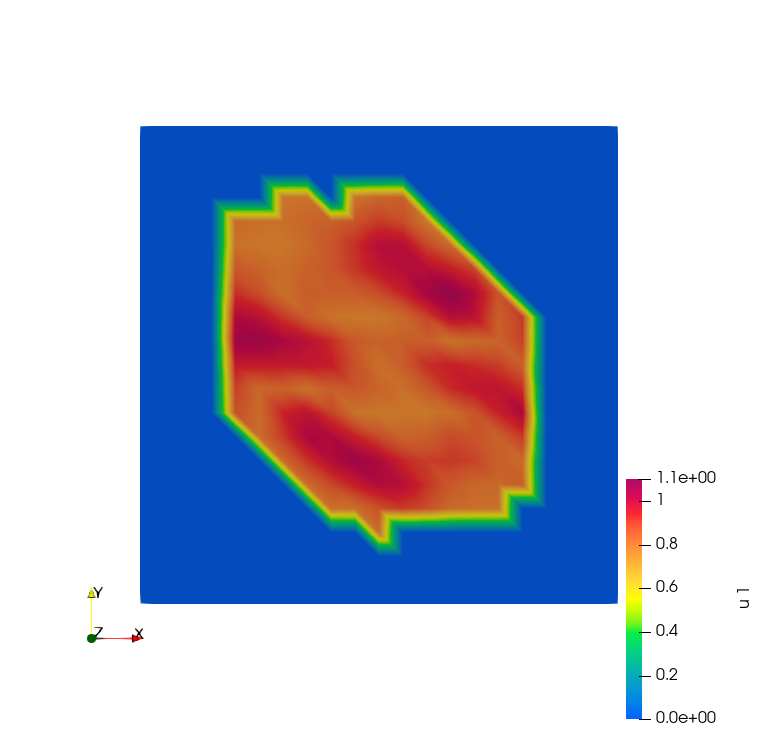}
      \\
      a) $\max \varepsilon_{r_h}/5  \approx 9$                                                                     & b)  $\max \sigma_{h}/5 \approx 1.1$ \\
      \includegraphics[trim = 1.5cm .2cm .8cm 3.5cm, scale=0.11, clip=]{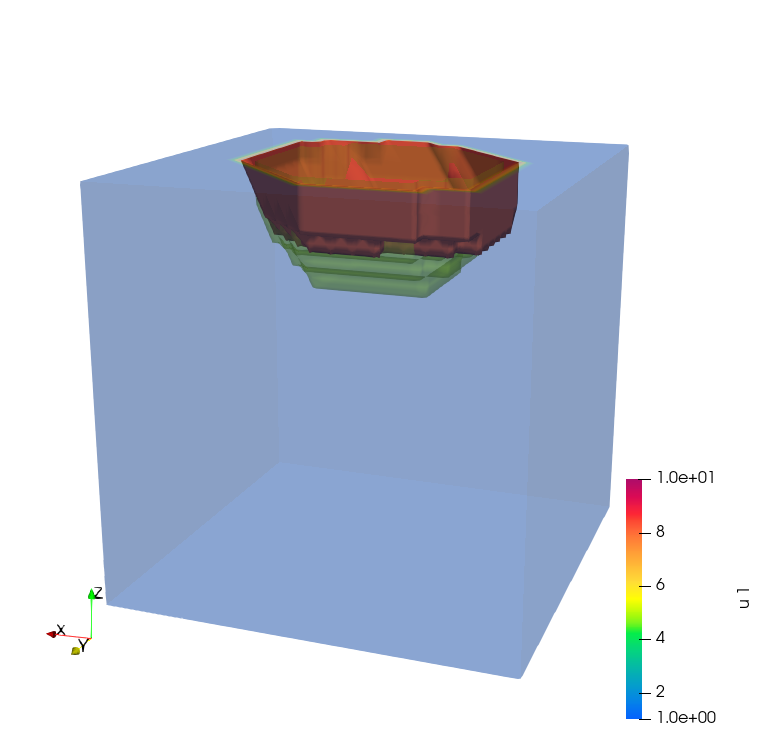}    &
      \includegraphics[trim = 2.0cm 0.2cm 2.81cm 3.0cm, scale=0.11, clip=]{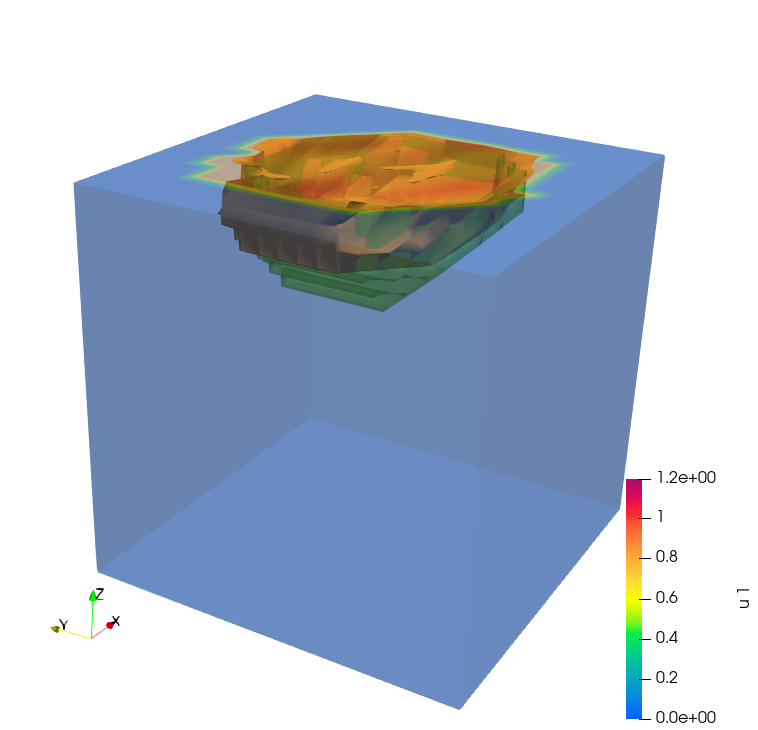}
      \\
      c) $\max \varepsilon_{r_h}/5  \approx 9$                                                                     & d)  $\max \sigma_{h}/5 \approx 1.1$ \\
      \includegraphics[trim = 11.0cm 6.2cm 11.5cm 6.2cm, scale=0.11, clip=]{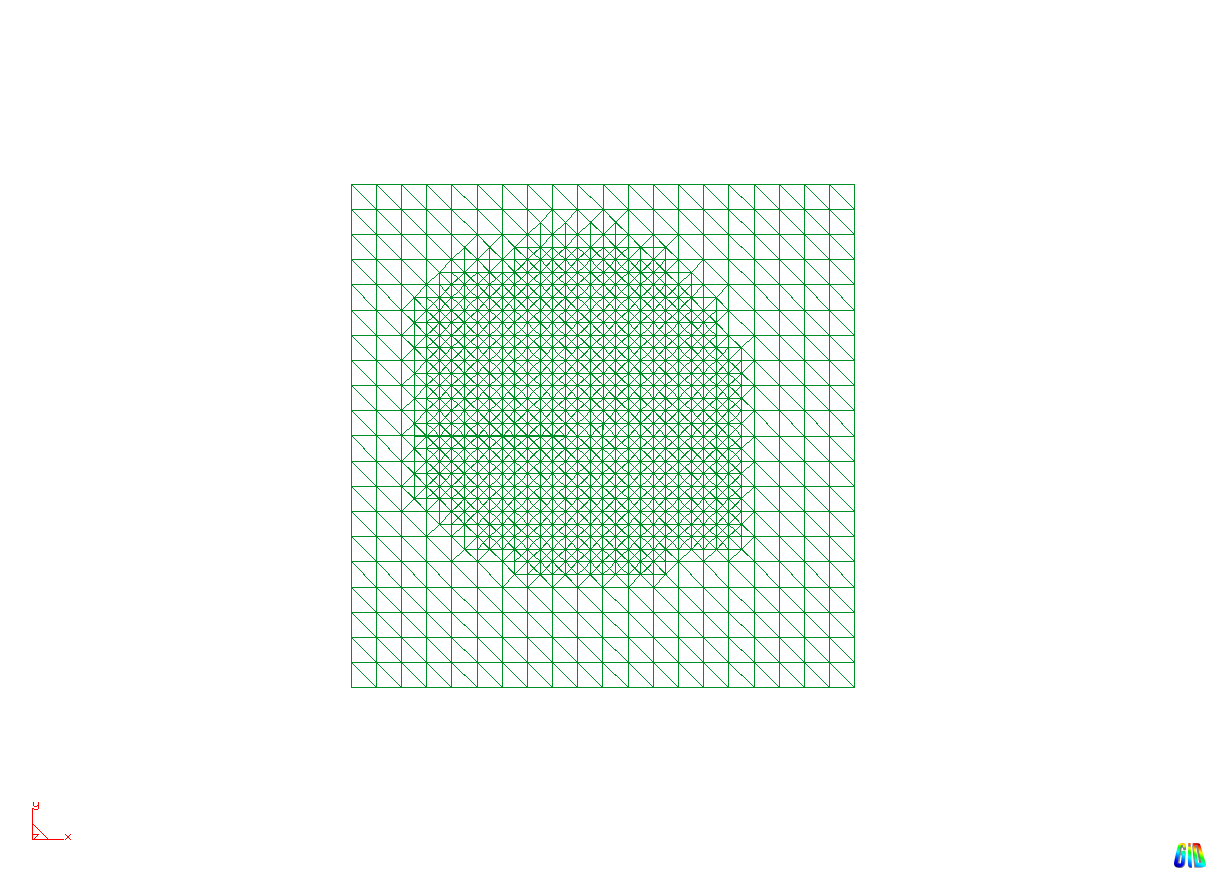}              &
      \includegraphics[trim = 11.0cm 6.2cm 11.5cm 6.2cm, scale=0.11, clip=]{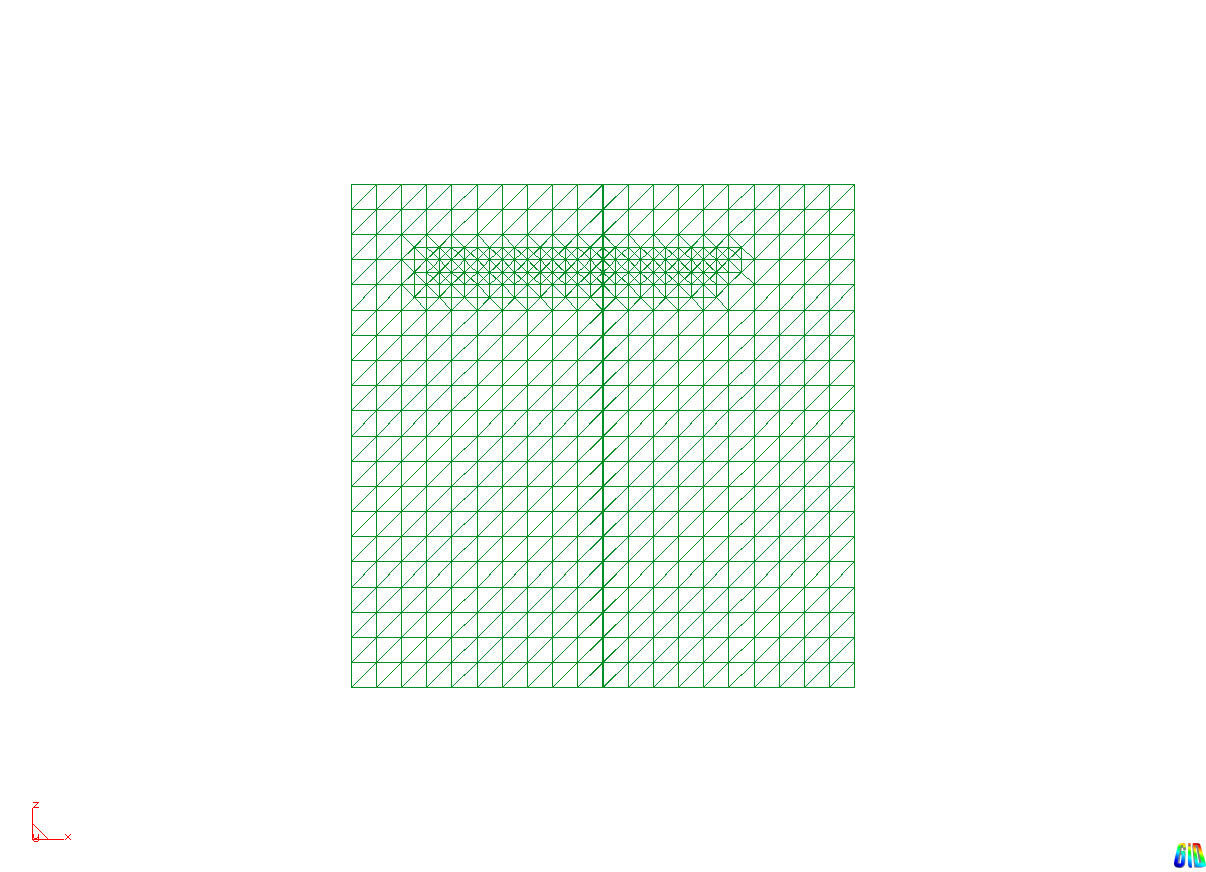}                                                    \\
      e)  $x_1 x_2$-view                                                                                           & f) $x_3 x_1$- view
    \end{tabular}
  \end{center}
  \caption{\small \textit{Performance of ACGA.
    a), c) The weighted reconstruction of $\varepsilon_{r_h}$ (outlined in red color) corresponding to
    the malign melanoma at stage 1 and
    obtained
    on six times locally refined finite element mesh ${K_{h_6}}$.
    b),d) The weighted reconstruction of $\sigma_{h}$ (outlined in red color) corresponding to
    the malign melanoma at stage 1 and
    obtained
    on five times locally refined finite element mesh ${K_{h_5}}$.
    c), d) Projections of the refined mesh ${K_{h_6}}$.
    Figures c), d) also present exact isosurface of melanoma (in yellow color) for comparison with reconstrcuted values of $\varepsilon_{r_h}/5, \sigma_{h}/5.$
    The noise level in the data for electric field is $\delta= 10\%$.}
  }
  \label{fig:ACGA}
\end{figure}

The goal of our numerical tests is to reconstruct weighted dielectric
permittivity $\varepsilon/5$ and conductivity $\sigma/5$ functions
presented in Figure \ref{fig:numex2}- a), b) using time-dependent
backscattered electric field $E =(E_1,E_2,E_3)$ at $\Gamma_{1}$.
 Figure \ref{fig:numex2} - c), d)
shows projections of locally refined finite element mesh for generation of
backscattered data.
Figure
\ref{fig:data} demonstrates the computed scattered electric field
$|E|$ of the model problem \eqref{model1} in the finite element
domain $\Omega_{\rm FEM} $ at different times using ADDFE/FDM method on the geometry
shown on Figure \ref{fig:meshes}.
Figures \ref{fig:CGA}, \ref{fig:ACGA} show reconstruction of the
weighted values of the relative dielectric permittivity and
conductivity functions obtained by CGA and ACGA algorithms,
correspondingly. We take as an initial guess $\varepsilon^0 = 1,
 \sigma^0 = 0$ for all points of the computational domain $\Omega$.
Comparing results of the reconstruction presented on Figure
\ref{fig:CGA} with reconstructions shown on Figure \ref{fig:ACGA}
we can conclude that the local adaptive mesh refinement, or ACGA
algorithm, significantly improve results of 3D reconstruction of
$\varepsilon$ and $\sigma$ and converges to the real depth of the
location of MM. We also observe that shape of the obtained
reconstructions is also correctly represented.

\section{Conclusions}

\label{sec:concl}

The work presents performance of CGA and ACGA algorithms for reconstruction of
dielectric properties of malign melanoma placed in the homogeneous domain.
Our computational tests show qualitative and quantitative
reconstruction of the relative dielectric permittivity function of
malign melanoma measured at 6 GHz.
Future work is related to the reconstruction of dielectric properties of malign melanoma placed inside skin with dielectric properties of skin corresponding to the real values of MM presented in the Table 1 and Figure \ref{fig:numex1}.

\section*{Acknowledgment}

The research of authors  is supported by the Swedish Research Council grant VR 2024-04459
 and STINT grant MG2023-9300.

\vspace{12pt}

\end{document}